\pdfoutput=1
\documentclass[a4paper,abstract]{scrartcl}

\usepackage[T1]{fontenc}
\usepackage[utf8]{inputenc}
\usepackage[english]{babel}
\usepackage{microtype}

\usepackage{graphicx}
\usepackage{amsmath}
\usepackage{amsfonts}
\usepackage{adjustbox}
\usepackage{mathtools}
\usepackage{pgfplots}
\pgfplotsset{compat=1.16}
\usepgfplotslibrary{groupplots}
\usepackage{catchfile}
\usepackage{booktabs}
\usepackage{algorithm}
\usepackage{algpseudocode}

\newcommand{\Span}{\operatorname{span}}
\newcommand{\im}{\operatorname{Im}}

\begin{document}

\title{Localized Reduced Basis Additive Schwarz Methods}
\author{
    Martin J. Gander\thanks{
        Section de Math\'ematiques, Universit\'e de Gen\`eve,
        2-4 rue du Li\`evre, CP 64, CH-1211, Gen\`eve, Suisse,
        \texttt{Martin.Gander@unige.ch}
    }
    \and
    Stephan Rave\thanks{
        Mathematics Münster, University of Münster, Einsteinstrasse 62, 48149 Münster, Germany,
        \texttt{Stephan.Rave@uni-muenster.de}
    }
}
\date{June 8, 2021}
\maketitle

\begin{abstract}
  \noindent Reduced basis methods build low-rank approximation spaces
  for the solution sets of parameterized PDEs
  by computing solutions of the given PDE for appropriately selected snapshot parameters.
  Localized reduced basis methods reduce the offline cost of computing these snapshot solutions
  by instead constructing a global space from spatially localized
  less expensive problems.
  In the case of online enrichment, these local problems are
  iteratively solved in regions of high residual and correspond to subdomain
  solves in domain decomposition methods.
  We show in this note that indeed there is a close relationship between
  online-enriched localized reduced basis and domain decomposition methods by introducing a
  Localized Reduced Basis Additive Schwarz method (LRBAS), which can be
  interpreted as a locally adaptive multi-preconditioning scheme for the
  CG method.
\end{abstract}

\section{Introduction}

Reduced basis (RB) methods~\cite{QuarteroniManzoniEtAl2016,HesthavenRozzaEtAl2016} are a family of model order reduction
schemes for parameterized PDEs, which can speed up the repeated
solution of such equations by orders of magnitude.
In the so-called \emph{offline} phase, RB methods construct a problem-adapted low-dimensional approximation space by
computing solutions of the PDE for selected \emph{snapshot} parameters using a given high-fidelity discretization of the PDE.
In the following \emph{online} phase, the PDE is solved for arbitrary new parameters by computing the (Petrov-)Galerkin
projection of its solution onto the precomputed reduced approximation space.
While RB methods have been proven successful in various applications, for very large problems the computation of the
solution snapshots in the offline phase may still be prohibitively expensive.
To mitigate this issue, localized RB methods~\cite{MADAY2002195, BuhrIapichinoEtAl2020} have been developed which construct the
global approximation space from spatially localized less expensive problems.
These local problems largely fall into two classes:

\emph{Training} procedures construct local approximation spaces without knowledge of the global problem by, e.g., solving
the equation on an enlarged subdomain with arbitrary boundary values and then restricting the solution to the domain of
interest, or by solving related eigenvalue problems.
As such, these training approaches have a strong connecting with numerical multiscale methods and the construction of
spectral coarse spaces in domain decomposition methods.

In this contribution, however, we will focus on the construction of local RB spaces via \emph{online enrichment}, where
these spaces are iteratively built by solving localized corrector problems for the residual of the current reduced
solution.
In particular we mention the use of online enrichment in context of the LRBMS~\cite{OhlbergerSchindler2015},
GMsFEM~\cite{ChungEfendievEtAl2015} and ArbiLoMod~\cite{BuhrEngwerEtAl2017} methods.
These enrichment schemes share strong similarities with Schwarz methods, and it is the main goal of this contribution to
shed some light on the connections between these methods.
We will do so by introducing a simple localized RB additive Schwarz (LRBAS) method which is phrased in the language of
the abstract Schwarz framework but incorporates the central ingredients of online adaptive localized RB methods.
In particular, we hope that LRBAS will help the analysis of localized RB methods from the perspective of Schwarz methods.
Following~\cite{BuhrEngwerEtAl2014}, we will consider arbitrary but localized changes of the problem instead of
parametric variations.
In Section~\ref{sec:mpcg} we will see that LRBAS can indeed be interpreted as a locally adaptive version of a multi-preconditioned CG
method.

Compared to Schwarz methods, a distinctive feature of LRBAS is that updates are only computed in high-residual regions,
which can lead to a significant reduction of the number of local updates and a concentration of the updates to a few
regions affected by the localized changes (cf.\ Section~\ref{sec:experiments}). 
This property might be exploited for the reduction of the overall power consumption and to balance the computational load among a
smaller amount of compute nodes, in particular in cloud environments, where additional computational resources can be
easily allocated and deallocated again.

\section{A Localized Reduced Basis Additive Schwarz Method}

Our goal is to efficiently solve a sequence, indexed by $k$, of linear systems
\begin{equation}\label{eq:fom}
    A^{(k)}x^{(k)} = f
\end{equation}
with $A^{(k)} \in \mathbb{R}^{n\times n}$ symmetric, positive definite and $x^{(k)}, f \in \mathbb{R}^n$, up to some fixed error
tolerance $\varepsilon$.
To this end, let $n\times n_i$ matrices $R_i^T$ of rank $n_i$ be given for $1 \leq i \leq I$
and $n\times n_0^{(k)}$ matrices $R_0^{(k)T}$ of rank $n_0^{(k)}$.
Typically, $R_1, \ldots, R_I$ will be the restriction matrices corresponding to a finite element basis associated with
an overlapping domain decomposition $\Omega_i$ of the computational domain $\Omega$, and the columns of $R_0^{(k)T}$ contain a
basis of a suitable coarse space for $A^{(k)}$.
In particular we assume that each $R_i$ is non-orthogonal to only a few neighboring spaces, i.e., there are a small
constant $C$ and index sets $\mathcal{O}_i \subset \{1,\ldots,I\}$ with $\#\mathcal{O}_i \leq C\cdot I$ such that
\begin{equation}\label{eq:locality}
    R_j\cdot R_i^T = 0_{n_j\times n_i} \qquad\text{whenever } j \notin \mathcal{O}_i.
\end{equation}
As usual, we define the local matrices
\begin{equation}\label{eq:local_operators}
    A_0^{(k)}:=R_0^{(k)}A^{(k)}R_0^{(k)T} \qquad\text{and}\qquad A_i^{(k)} := R_iA^{(k)}R_i^T.
\end{equation}
We are interested in the case where $A^{(k+1)}$ is obtained from $A^{(k)}$ by an arbitrary but local modification in the
sense that
\begin{equation}\label{eq:local_change}
    A_i^{(k+1)} = A_i^{(k)} \qquad\text{for } i\notin \mathcal{C}^{(k+1)},
\end{equation}
where the sets $\mathcal{C}^{(k)}$ contain the indices of the spaces affected by the change, generally assuming
that $\# \mathcal{C}^{(k)} \ll I$.

Over the course of the computation of the solutions $x^{(k)}$ we will build local low-dimensional reduced bases 
$\tilde{R}_i^{(k,l)T} \in \mathbb{R}^{n_i\times N_{i}^{(k,l)}}$ for $i \geq 1$ such that there are local coefficients
$\tilde{x}_i^{(k,l)} \in \mathbb{R}^{N_i^{(k,l)}}$ and $\tilde{x}_0^{(k,l)} \in \mathbb{R}^{n_0^{(k)}}$ such that
\begin{equation}\label{eq:decomposition}
    \tilde{x}^{(k,l)} := R_0^{(k)T} \tilde{x}^{(k,l)}_0 + \sum_{i=1}^{I} R_i^T \tilde{R}_i^{(k,l)T} \tilde{x}^{(k,l)}_i
\end{equation}
is a good approximation of $x^{k}$ for sufficiently large $l$.
We obtain such an approximation via Galerkin projection onto the global reduced basis space spanned by the images of
$R_0^{(k)T}$ and all $R_i^T\tilde{R}_i^{(k,l)T}$, i.e., $\tilde{x}^{(k,l)}$ is determined by the
$(n_0^{(k)} + \sum_{i=1}^I N_i^{(k,l)})$-dimensional linear system
\begin{equation}\label{eq:rom}
    \begin{aligned}
        R^{(k)}_0 A^{(k)}\tilde{x}^{(k,l)} &= R^{(k)}_0 f, \\
        \tilde{R}_i^{(k,l)}R_i A^{(k)} \tilde{x}^{(k,l)} &= \tilde{R}_i^{(k,l)}R_i f, &\quad 1\leq i \leq I.
    \end{aligned}
\end{equation}
Thanks to the locality~\eqref{eq:locality} of the space decomposition, the matrix of the system \eqref{eq:rom}
has a block structure allowing us to efficiently assemble and solve it. 

To build the local reduced bases $\tilde{R}_i^{(k,l)T}$ we use an iterative enrichment procedure where the basis is
extended with local Schwarz corrections $y_i^{(k,l)} \in \mathbb{R}^{n_i}$ for the current residual,
\begin{equation}\label{eq:enrichment}
    A^{(k)}_i y_i^{(k,l)} = r_i^{(k,l)} := R_i (f - A^{(k)} \tilde{x}^{(k)}).
\end{equation}
In view of~\eqref{eq:local_change}, the corrections are only computed in subdomains $i$ with large residual
norm $\|r_i^{(k,l)}\|$.
In particular, for
finite-element discretizations of elliptic PDEs without
high-conductivity channels, we expect that with increasing $k$ the number of enriched
bases will be of the same order as the cardinality of $\mathcal{C}^{k+1}$.
The exact definition of the enrichment scheme is given in Algorithm~\ref{alg:lrbas}.
There are various possibilities to choose the criterion for the localized enrichment in line 9 of Algorithm~\ref{alg:lrbas}.
In this work we simply select those reduced spaces for enrichment for which the quotient between the norm of
the local residual and the norm of the global residual is larger than a fixed constant that scales with the number of
the subdomains.

\begin{algorithm}[t]
    \caption{\underline{L}ocalized \underline{R}educed \underline{B}asis \underline{A}dditive \underline{S}chwarz method
    (LRBAS)}\label{alg:lrbas}
    \begin{algorithmic}[1]
        \Procedure{LRBAS}{$A^{(k)}$, $f$, $R_0^T$, $R_i^T$, $\varepsilon$, $\varepsilon_{\mathrm{loc}}$}
        \State{$\tilde{R}_i^{(1,1)T} \gets 0_{n_i\times 0}, \quad 1\leq i \leq I$}
        \Comment{initialize local bases}
        \For{$k \gets 1, \ldots, \infty$}
        \State{$\tilde{x}^{(k,1)}, \tilde{x}_i^{(k,1)} \gets \text{solutions of \eqref{eq:decomposition}, \eqref{eq:rom}}$}
        \Comment{initial solution}
        \State{$r^{(k,1)} \gets f - A^{(k)} \tilde{x}^{(k,1)}$} \Comment{initial residual}
        \State{$l \gets 1$}
        \While{$\|r^{(k,l)}\|\,/\, \|f\| > \varepsilon$}\Comment{loop until converged}
        \For{$i \gets 1, \ldots I$}\Comment{enrichment procedure}
        \If{$\|R_i r^{(k,l)}\|^2 > \varepsilon_{\mathrm{loc}} \cdot I^{-1} \cdot \|r^{(k,l)}\|^2$}
        \State{$y_i^{(k,l)} \gets \text{solution of~\eqref{eq:enrichment}}$}
        \State{$\tilde{R}_i^{(k,l+1)T} \gets \bigl[\tilde{R}_i^{(k,l)T}\ \, y_i^{(k,l)} \bigr]$}
        \Else
        \State{$\tilde{R}_i^{(k,l+1)} \gets \tilde{R}_i^{(k,l)}$}
        \EndIf{}
        \EndFor
        \State{$\tilde{x}^{(k,l+1)}, \tilde{x}_i^{(k,l+1)} \gets \text{solutions of \eqref{eq:decomposition}, \eqref{eq:rom}}$}
        \Comment{update solution}
        \State{$r^{(k,l+1)} \gets f - A^{(k)} \tilde{x}^{(k,l+1)}$}
        \Comment{update residual}
        \State{$l \gets l+1$}
        \EndWhile{}
        \For{$i \gets 1, \ldots I$}\Comment{update bases for next problem}
        \If{$\tilde{R}_i^{(k,l)T} \neq \tilde{R}_i^{(k,1)T}$}\Comment{basis enriched at least once?}
        \State{$\tilde{R}_i^{(k+1,1)T} \gets \bigl[\tilde{R}_i^{(k,1)T}\ \, \tilde{R}_i^{(k,l)T}\tilde{x}_i^{(k,l)}\bigr]$}\Comment{only keep local solution in basis}
        \Else
        \State{$\tilde{R}_i^{(k+1,1)T} \gets \tilde{R}_i^{(k,1)T}$}
        \EndIf
        \EndFor
        \EndFor
        \EndProcedure
    \end{algorithmic}
\end{algorithm}

Note that an important property of localized enrichment is that after an enrichment step only those blocks $(i,j)$ of the matrix
corresponding to~\eqref{eq:rom} have to be updated for which either $\tilde{R}_i^{(k,l)}$ or $\tilde{R}_j^{(k,l)}$ have been
enriched.
Using reduced basis techniques~\cite{BuhrEngwerEtAl2014} it is further possible to evaluate the residual norms $\|R_i r^{(k,l)}\|$ and
$\|r^{(k,l)}\|$ using only reduced quantities, which again only have to be updated for local bases $\tilde{R}_i^{(k,l)T}$
affected by the enrichment.
Thus, in a distributed computing environment only the main compute node solving~\eqref{eq:rom} and those nodes
associated with the enriched bases have to perform any operations, while the other compute node lay at rest.

We remark that several extensions to the LRBAS method are possible.
In particular, we assumed for simplicity that all matrices $A^{(k)}$ are of the same dimension.
This, for instance, is the case when coefficient functions of the PDE underlying~\eqref{eq:fom} are modified, but
the computational mesh remains unchanged.
However, also local geometry changes that lead to remeshing can be handled by resetting all local bases that are
supported on the changed geometry.
In this context we note that, as another simplification, in the definition of LRBAS we have chosen to keep all basis
vectors when transitioning from $A^{(k)}$ to $A^{(k+1)}$, including bases $\tilde{R}_i^{(k,l)T}$ affected by the change,
even though these retained bases will generally not contribute to the convergence of the
scheme.
Finally, in many applications, a local or global parametric variation of $A^{(k)}$, e.g.~the change of some material
parameters, in addition to the considered non-parametric modifications may be of interest.
In such cases, parametric model order reduction techniques such as greedy basis generation algorithms or offline/online
decomposition of the reduced order system~\eqref{eq:rom} can be incorporated into the scheme.
In particular we refer to~\cite{BuhrEngwerEtAl2017} where both additional parameterization of $A^{(k)}$ as well as the
reinitialization of the local bases after non-parametric changes from $A^{(k)}$ to $A^{(k+1)}$ are discussed.

\subsection{LRBAS as an additive-Schwarz multi-preconditioned CG method}\label{sec:mpcg}
\newcommand{\Mkinv}{{\big(M^{(k)}\big)^{-1}}}
\newcommand{\SearchSpace}[1]{\mathcal{S}^{(k,l)}_{\text{#1}}}

Consider the solution of the systems~\eqref{eq:fom} with the preconditioned conjugate gradient (PCG) algorithm, where we
choose as preconditioner the additive Schwarz operator
$\Mkinv := R_0^{(k)T}\big(A_0^{(k)}\big)^{-1}R_0^{(k)} + \sum_{i=1}^{I} R_i^T\big(A_i^{(k)}\big)^{-1}R_i$.
Let $x^{(k,l)}_\text{pcg}$ denote the $l$-th iterate of the PCG algorithm, starting with $x^{(k,0)} = 0$ as the initial
guess.
Then it is well known that $x^{(k,l)}_\text{pcg}$ lies in the search space $\SearchSpace{pcg}$ given by the Krylov space
$\mathcal{K}^{l}\left(\Mkinv A^{(k)}, \Mkinv f\right)$
and that the error $x^{(k)} - x^{(k,l)}_\text{pcg}$ is $A^{(k)}$-orthogonal to this space.
Denoting by $r^{(k,l)}_\text{pcg} := f - A^{(k)}x^{(k,l)}_{\text{pcg}}$ the $l$-th residual, one readily checks that
$\SearchSpace{pcg}$ is equivalently given by
\begin{equation}
    \SearchSpace{pcg} :=
    \Span \left\{\Mkinv r^{(k,0)}_{\text{pcg}},\; \dots,\; \Mkinv r^{(k,l-1)}_{\text{pcg}}\right\},
\end{equation}
i.e., in each iteration the search space is extended by the vector obtained from the application of the preconditioner
to the current residual.
The idea of multi-preconditioning~\cite{BridsonGreif2006} is to enlarge this search space by including each local preconditioner
$\big(A_i^{(k)}\big)^{-1}$ application into the search space individually, leading to
\begin{equation}
    \begin{multlined}
    \SearchSpace{mpcg} :=
        \Span 
        \Big(\left\{R_0^{(k)T}\big(A^{(k)}_0\big)^{-1}R_0^{(k)} r^{(k,t)}_{\text{mpcg}}\;\middle|\; 0 \leq t \leq l-1
        \right\}\\
        \qquad\qquad\qquad\qquad\cup
        \left\{R_i^T\big(A^{(k)}_i\big)^{-1}R_i r^{(k,t)}_{\text{mpcg}}\;\middle|\; 1\leq i \leq I,\, 0 \leq t \leq l-1 \right\}\Big),
    \end{multlined}
\end{equation}
with $r^{(k,l)}_{\text{mpcg}}$ denoting the multi-preconditioned CG residuals.
Conversely, we easily see from~\eqref{eq:rom} and~\eqref{eq:enrichment} that for $\varepsilon_{\text{loc}} = 0$ the
LRBAS iterates $\tilde{x}^{(k,l)}$ lie within the search space
\begin{equation}
    \begin{multlined}
        \SearchSpace{lrbas,0} :=
        \im\Bigl(\Bigl[R_0^{(k)T}\ \, R_1^{T}\tilde{R}_1^{(k,1)}\ \, \dots\ \, R_I^{T}\tilde{R}_I^{(k,1)}\Bigr]\Bigr)\\
        \qquad\qquad\qquad+ \Span \left\{R_i^{T}\big(A^{(k)}_i\big)^{-1}R_i r^{(k,t)}\;\middle|\; 1\leq i \leq I,\, 1 \leq t \leq l-1 \right\},
    \end{multlined}
\end{equation}
and that the error $x^{(k)} - \tilde{x}^{(k,l)}$ is $A^{(k)}$-orthogonal to this space.
Hence, LRBAS with $\varepsilon_{\text{loc}} = 0$ can be seen as a projected multi-preconditioned CG method for
solving~\eqref{eq:rom}, where the projection space is given by the span of the coarse space and the initial local
reduced bases and where the new solution iterate $\tilde{x}^{(k,l)}$ is obtained by direct solution of the reduced
system~\eqref{eq:rom} instead of an incremental update in order to preserve the locality of the reduced bases.

For $\varepsilon_{\text{loc}} > 0$ we arrive at an adaptive version of multi-preconditioning similar to~\cite{Spillane2016}.
However, in contrast to~\cite{Spillane2016} where either all local search directions or their global sum are added to
the search space, LRBAS is locally adaptive in the sense that only those local search directions are computed and included where a
large local residual has to be corrected.

\section{Numerical Experiment}\label{sec:experiments}

We consider the test case from~\cite{BuhrEngwerEtAl2017} and solve a sequence of five elliptic problems
\begin{equation}\label{eq:testcase}
    \begin{aligned}
        \nabla \cdot \big(- \sigma^{(k)}(x,y) \nabla u^{(k)}(x,y)\big) &= 0, & x,y&\in (0,1),\\
        u^{(k)}(0, y) &= 1, & y &\in (0,1),\\
        u^{(k)}(1, y) &= -1, & y &\in (0,1),\\
        -\sigma^{(k)}(x,y) \nabla u^{(k)}(x,y)\cdot \mathbf{n}(x,y) &= 0, & x&\in (0,1),\ y\in \{0,1\},
    \end{aligned}
\end{equation}
where the coefficient $\sigma^{(k)}(x)$ is given as in Fig.~\ref{fig:diffusivity}.
The problem is discretized using bilinear finite elements over a uniform $200 \times 200$ mesh.
The resulting solutions are visualized in Fig.~\ref{fig:solutions}.
We decompose the computational domain uniformly into $10 \times 10$ subdomains with an overlap of 4 mesh
elements.
For $R^{(k)}_0$ we choose GenEO~\cite{SpillaneDoleanEtAl2013} basis functions with eigenvalues below 0.5, yielding
between two and five functions per subdomain.
When connecting or disconnecting the high-conductivity channels, we expect enrichment to be required along
the subdomains adjacent to the channels, whereas the other subdomains should be largely unaffected by the local change.

\let\pgfimageWithoutPath\pgfimage 
\renewcommand{\pgfimage}[2][]{\pgfimageWithoutPath[#1]{plots/#2}}
\pgfplotsset{width=4.5cm, height=4.5cm,%
             every tick label/.append style={font=\scriptsize}}
\begin{figure}[t]
\hspace{0.7cm}
\adjustbox{valign=c}{%
\begin{tikzpicture}[baseline]
\begin{axis}[
colorbar,
colorbar style={ylabel={}},
colormap/viridis,
point meta max=100001,
point meta min=1,
tick align=outside,
tick pos=left,
x grid style={white!69.0196078431373!black},
xmin=0, xmax=1,
y grid style={white!69.0196078431373!black},
ymin=0, ymax=1,
ytick style={color=black}
]
\addplot graphics [includegraphics cmd=\pgfimage,xmin=0, xmax=1, ymin=0, ymax=1] {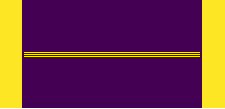};
\node (A) [coordinate] at (axis cs:0.105, 0.52) {};
\node (B) [coordinate] at (axis cs:0.105, 0.5) {};
\node (C) [coordinate] at (axis cs:0.105, 0.48) {};
\node (D) [coordinate] at (axis cs:0.892, 0.52) {};
\node (E) [coordinate] at (axis cs:0.892, 0.5) {};
\node (F) [coordinate] at (axis cs:0.892, 0.48) {};
\node [inner sep=0, red] (Atext) at (axis cs:0.05, 0.58) {1};
\node [inner sep=0, red] (Btext) at (axis cs:0.05, 0.5)  {2};
\node [inner sep=0, red] (Ctext) at (axis cs:0.05, 0.42) {3};
\node [inner sep=0, red] (Dtext) at (axis cs:0.95, 0.58) {4};
\node [inner sep=0, red] (Etext) at (axis cs:0.95, 0.5)  {5};
\node [inner sep=0, red] (Ftext) at (axis cs:0.95, 0.42) {6};
\draw [red, semithick] (Atext) -- (A);
\draw [red, semithick] (Btext) -- (B);
\draw [red, semithick] (Ctext) -- (C);
\draw [red, semithick] (Dtext) -- (D);
\draw [red, semithick] (Etext) -- (E);
\draw [red, semithick] (Ftext) -- (F);
\end{axis}
\end{tikzpicture}%
}
\hfill
\adjustbox{valign=c}{%
\setlength{\tabcolsep}{4pt}%
\newcommand{\C}{\checkmark}
\begin{tabular}{rccccc}
\toprule
  $k=$     &  1  &  2  &  3  &  4  &  5  \\\midrule
 ch. 1     &     &     &     &  \C &  \C \\
 ch. 2     &  \C &     &     &     &     \\
 ch. 3     &     &     &     &     &     \\
 ch. 4     &     &     &     &     &     \\
 ch. 5     &  \C &  \C &     &     &  \C \\
 ch. 6     &     &     &     &     &     \\\bottomrule
\end{tabular}%
\hspace{1.0cm}
}
\caption{Definition of the coefficient functions $\sigma^{(k)}$ for the numerical test case~\eqref{eq:testcase};
left: function $\sigma^{(0)}$, taking the values $10^5 + 1$ inside the high-conductivity regions
and $1$ elsewhere; right: $\sigma^{(k)}$ is obtained from $\sigma^{(0)}$ by connecting the three channels to the boundary
regions at the marked locations.}\label{fig:diffusivity}
\end{figure}
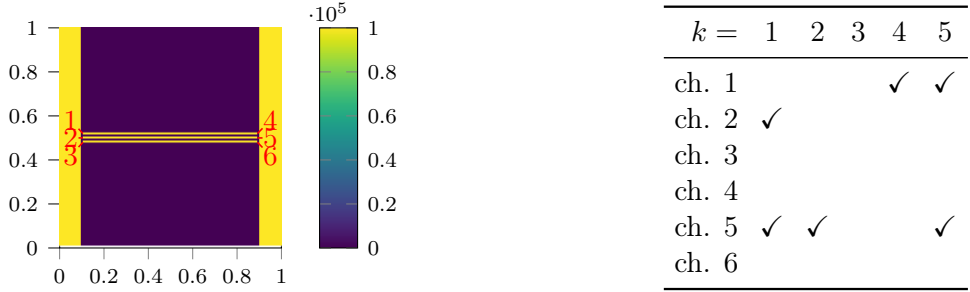%
\pgfplotsset{width=4.35cm, height=4.35cm,%
             every tick label/.append style={font=\scriptsize}}
  \begin{figure}[t]
    \begin{center}
\begin{tikzpicture}

\begin{groupplot}[group style={group size=3 by 2}]
\nextgroupplot[
tick align=outside,
tick pos=left,
x grid style={white!69.0196078431373!black},
xmin=0, xmax=1,
xtick style={color=black},
y grid style={white!69.0196078431373!black},
ymin=0, ymax=1,
ytick style={color=black}
]
\addplot graphics [includegraphics cmd=\pgfimage,xmin=0, xmax=1, ymin=0, ymax=1] {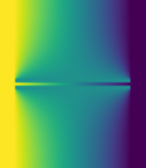};

\nextgroupplot[
tick align=outside,
tick pos=left,
x grid style={white!69.0196078431373!black},
xmin=0, xmax=1,
xtick style={color=black},
y grid style={white!69.0196078431373!black},
ymin=0, ymax=1,
ytick style={color=black}
]
\addplot graphics [includegraphics cmd=\pgfimage,xmin=0, xmax=1, ymin=0, ymax=1] {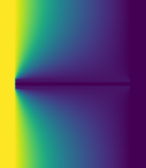};

\nextgroupplot[
tick align=outside,
tick pos=left,
x grid style={white!69.0196078431373!black},
xmin=0, xmax=1,
xtick style={color=black},
y grid style={white!69.0196078431373!black},
ymin=0, ymax=1,
ytick style={color=black}
]
\addplot graphics [includegraphics cmd=\pgfimage,xmin=0, xmax=1, ymin=0, ymax=1] {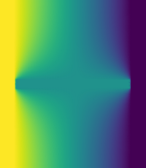};

\nextgroupplot[
tick align=outside,
tick pos=left,
x grid style={white!69.0196078431373!black},
xmin=0, xmax=1,
xtick style={color=black},
y grid style={white!69.0196078431373!black},
ymin=0, ymax=1,
ytick style={color=black}
]
\addplot graphics [includegraphics cmd=\pgfimage,xmin=0, xmax=1, ymin=0, ymax=1] {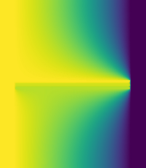};

\nextgroupplot[
colorbar,
colorbar style={ylabel={}},
colormap/viridis,
point meta max=1,
point meta min=-1,
tick align=outside,
tick pos=left,
x grid style={white!69.0196078431373!black},
xmin=0, xmax=1,
xtick style={color=black},
y grid style={white!69.0196078431373!black},
ymin=0, ymax=1,
ytick style={color=black}
]
\addplot graphics [includegraphics cmd=\pgfimage,xmin=0, xmax=1, ymin=0, ymax=1] {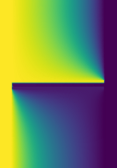};
\end{groupplot}

\end{tikzpicture}
    \end{center}
    \caption{Solutions of the test problem~\eqref{eq:testcase} for $k=1,2,3$ (top row) and $k=4,5$ (bottom row).}\label{fig:solutions}
\end{figure}
{
\newcommand{\PCGit}{107}
\newcommand{\PCGls}{10700}
\newcommand{\PCGGuessit}{63}
\newcommand{\PCGGuessls}{6300}
\newcommand{\MPCGit}{33}
\newcommand{\MPCGls}{3300}
\newcommand{\MPCGKeepit}{28}
\newcommand{\MPCGKeepls}{2800}
\newcommand{\AMPCGit}{39}
\newcommand{\AMPCGls}{1386}
\newcommand{\AMPCGKeepit}{34}
\newcommand{\AMPCGKeepls}{1335}
\newcommand{\VNT}{}
\begin{table}
    \begin{center}
        \setlength{\tabcolsep}{6pt}
        \begin{tabular}{lcc}
        \toprule
                                                              & iterations     & local enrichments~\eqref{eq:enrichment} \\\midrule
        PCG                                                   & 107{}       & 10700{}        \\
        PCG + LRB solution as initial value                    & 63{}  & 6300{}   \\
        LRBAS ($\varepsilon_\text{loc} = 0$)                  & 33{}      & 3300{}       \\
        LRBAS ($\varepsilon_\text{loc} = 0.25$)                & 39{}     & 1386{}      \\
        LRBAS ($\varepsilon_\text{loc} = 0$, $\tilde{R}_i^{(k+1,1)} := \tilde{R}_i^{(k,l)}$) & 28{}  & 2800{}   \\
        LRBAS ($\varepsilon_\text{loc} = 0.25$, $\tilde{R}_i^{(k+1,1)} := \tilde{R}_i^{(k,l)}$) & 34{} & 1335{}  \\\bottomrule
        \end{tabular}
    \end{center}
    \caption{Total number of iterations and local Schwarz corrections~\eqref{eq:enrichment} required to reach a relative
    error tolerance $\varepsilon = 10^{-6}$ for the test problem~\eqref{eq:testcase}.}\label{tab:iterations}
\end{table}}
In Table~\ref{tab:iterations} we compare the total number of iterations for all five problems and the total
number of Schwarz corrections~\eqref{eq:enrichment} required to reach a relative error tolerance of $\varepsilon = 10^{-6}$
for the following solution strategies: 1.~the additive Schwarz preconditioned CG method with zero initial guess or with a localized RB
solution as initial guess, where the localized basis is obtained from the linear span of previous solutions $x^{(k)}$
decomposed using the GenEO partition of unity; 2.~LRBAS with and without local adaptivity ($\varepsilon_{\text{loc}} = 0.25$
or $0$); 3.~a version of LRBAS where the entire bases $\tilde{R}_i^{(k,l)T}$ are preserved when transitioning to
$k+1$ instead of only the final solution $\tilde{R}_i^{(k,l)T}\tilde{x}_i^{(k,l)}$.
As we see, LRBAS with locally adaptive enrichment significantly outperforms the PCG method with or without initial
guess, both regarding the number of required iterations as well as the number of Schwarz corrections. 
Compared to non-adaptive multi-preconditioning, i.e.\ LRBAS with $\varepsilon_\text{loc} = 0$, the number of
local corrections is more than halved at the expense of a slightly increased number of iterations.
Keeping all of $\tilde{R}_i^{(k,l)}$ improves the convergence of the method only slightly.
Finally, in Fig.~\ref{fig:basis_sizes}
\begin{figure}[t]
    \begin{center}
\begin{tikzpicture}

\begin{groupplot}[group style={group size=3 by 2}]
\nextgroupplot[
tick align=outside,
tick pos=left,
x grid style={white!69.0196078431373!black},
xmin=0, xmax=1,
xtick style={color=black},
y grid style={white!69.0196078431373!black},
ymin=0, ymax=1,
ytick style={color=black}
]
\addplot graphics [includegraphics cmd=\pgfimage,xmin=0, xmax=1, ymin=0, ymax=1] {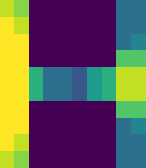};

\nextgroupplot[
tick align=outside,
tick pos=left,
x grid style={white!69.0196078431373!black},
xmin=0, xmax=1,
xtick style={color=black},
y grid style={white!69.0196078431373!black},
ymin=0, ymax=1,
ytick style={color=black}
]
\addplot graphics [includegraphics cmd=\pgfimage,xmin=0, xmax=1, ymin=0, ymax=1] {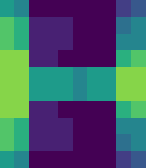};

\nextgroupplot[
tick align=outside,
tick pos=left,
x grid style={white!69.0196078431373!black},
xmin=0, xmax=1,
xtick style={color=black},
y grid style={white!69.0196078431373!black},
ymin=0, ymax=1,
ytick style={color=black}
]
\addplot graphics [includegraphics cmd=\pgfimage,xmin=0, xmax=1, ymin=0, ymax=1] {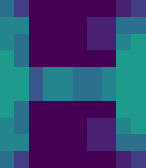};

\nextgroupplot[
tick align=outside,
tick pos=left,
x grid style={white!69.0196078431373!black},
xmin=0, xmax=1,
xtick style={color=black},
y grid style={white!69.0196078431373!black},
ymin=0, ymax=1,
ytick style={color=black}
]
\addplot graphics [includegraphics cmd=\pgfimage,xmin=0, xmax=1, ymin=0, ymax=1] {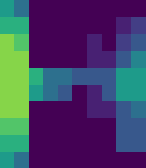};

\nextgroupplot[
colorbar,
colorbar style={ylabel={}},
colormap/viridis,
point meta max=11,
point meta min=0,
tick align=outside,
tick pos=left,
x grid style={white!69.0196078431373!black},
xmin=0, xmax=1,
xtick style={color=black},
y grid style={white!69.0196078431373!black},
ymin=0, ymax=1,
ytick style={color=black}
]
\addplot graphics [includegraphics cmd=\pgfimage,xmin=0, xmax=1, ymin=0, ymax=1] {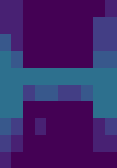};
\end{groupplot}

\end{tikzpicture}
    \end{center}
    \caption{Number of local Schwarz corrections~\eqref{eq:enrichment} required by the LRBAS method with
        $\varepsilon_\text{loc} = 0.25$ to solve the five test problems~\eqref{eq:testcase} up to a relative error
    tolerance of $\varepsilon = 10^{-6}$.}\label{fig:basis_sizes}
\end{figure}
we depict the number of required Schwarz corrections per subdomain for each $k$.
We observe a good localization of the computational work among the subdomains most affected by the local changes.

\paragraph{Acknowledgement}
    Funded by the Deutsche Forschungsgemeinschaft (DFG, German Research Foundation) under Germany's Excellence Strategy EXC 2044 –390685587, Mathematics Münster: Dynamics–Geometry–Structure.

\end{document}